# ON NORMAL APPROXIMATIONS TO *U*-STATISTICS

By Vidmantas Bentkus, Bing-Yi Jing[1] and Wang Zhou[2]

*Vilnius Institute of Mathematics and Informatics, Hong Kong University of Science and Technology and National University of Singapore*

Let $X_1, \ldots, X_n$ be i.i.d. random observations. Let $\mathbb{S} = \mathbb{L} + \mathbb{T}$ be a *U*-statistic of order $k \geq 2$ where $\mathbb{L}$ is a linear statistic having asymptotic normal distribution, and $\mathbb{T}$ is a stochastically smaller statistic. We show that the rate of convergence to normality for $\mathbb{S}$ can be simply expressed as the rate of convergence to normality for the linear part $\mathbb{L}$ plus a correction term, $(\operatorname{var} \mathbb{T}) \ln^2(\operatorname{var} \mathbb{T})$, under the condition $\mathbb{E}\mathbb{T}^2 < \infty$. An optimal bound without this log factor is obtained under a lower moment assumption $\mathbb{E}|\mathbb{T}|^\alpha < \infty$ for $\alpha < 2$. Some other related results are also obtained in the paper. Our results extend, refine and yield a number of related-known results in the literature.

**1. Introduction.** There has been a vast literature related to normal approximations and to the rates of convergence to normality for *U*-statistics of order $k \geq 2$. Undoubtedly, the case for $k = 2$ has been most studied and is best understood, so we will start our discussion with this case as well. Let $X, X_1, X_2, \ldots, X_n$, $n \geq 2$, be independent and identically distributed (i.i.d.) random variables (r.v.'s). Define a *U*-statistic of order 2 by

$$(1.1) \qquad U_n = \frac{2}{n(n-1)} \sum_{1 \leq i < j \leq n} h(X_i, X_j),$$

where the kernel $h(x, y)$ is a real-valued Borel measurable function, symmetric in its arguments with $Eh(X_1, X_2) = \theta$. By the Hoeffding decomposition, we have

$$\frac{\sqrt{n}(U_n - \theta)}{2\sigma_g} = \frac{1}{\sqrt{n}\sigma_g} \sum_{i=1}^n g(X_i) + \frac{1}{\sqrt{n}(n-1)\sigma_g} \sum_{1 \leq i < j \leq n} \eta(X_i, X_j)$$

Received December 2008; revised March 2009.
[1]Supported in part by Hong Kong RGC Grants HKUST6011/07P and HKUST6015/08P.
[2]Supported in part by Grant R-155-000-076-112 at National University of Singapore.

(1.2)
$$\stackrel{\text{def}}{=} \mathbb{L} + \mathbb{T},$$

where $\mathbb{L}$ and $\mathbb{T}$ are the linear and quadratic terms, respectively, and

$$g(x) = E[h(x, X_2)] - \theta, \qquad \sigma_g^2 = \text{Var}[g(X_1)],$$
$$\eta(x, y) = h(x, y) - g(x) - g(y) + \theta.$$

Throughout the paper, it is assumed that $\sigma_g^2 \in (0, \infty)$.

For any generic r.v. $Y$, denote its distribution function (d.f.) by $F_Y(x) = \mathbb{P}(Y \leq x)$. Then

$$F_{\mathbb{L}+\mathbb{T}}(x) = \mathbb{P}\left(\frac{\sqrt{n}(U_n - \theta)}{2\sigma_g} \leq x\right).$$

It is well known that $F_{\mathbb{L}+\mathbb{T}}(x)$ converges to the standard normal d.f., denoted by $\Phi(x)$, provided $Eh^2(X_1, X_2) < \infty$ (see [11]). In fact, this moment condition can be reduced to $Eg^2(X_1) < \infty$ and $E|h(X_1, X_2)|^{4/3} < \infty$ (see Remark 4.2.4 of [15], page 131). There also has been much work on the convergence rates to normality for $U$-statistics of order 2. However, the sharpest Berry–Esseen bound of order $O(n^{-1/2})$ comes from [14] and [9], who establish the following ideal bound:

(1.3)
$$\|F_{\mathbb{L}+\mathbb{T}} - \Phi\| \ll \frac{1}{\sqrt{n}}\left(\frac{E|g_1|^3}{\sigma_g^3} + \frac{E|\eta_{12}|^{5/3}}{\sigma_g^{5/3}}\right),$$

where the following notation has been used:

- $g_i = g(X_i)$ and $\eta_{ij} = \eta(X_i, X_j)$.
- For any function $f: \mathbb{R} \to \mathbb{R}$, define $\|f\| = \sup_{x \in \mathbb{R}} |f(x)|$.
- By $a \ll b$ for $a, b \geq 0$ we mean $a \leq cb$ for some positive constant $c$ not depending on the underlying distribution function.

Indeed, Bentkus, Götze and Zitikis [4] showed that the moment conditions in (1.3), $E|g_1|^3 < \infty$ and $E|\eta_{12}|^{5/3} < \infty$, are the weakest possible in the Berry–Esseen bounds of order $O(n^{-1/2})$ for $U$-statistics of order 2.

Our main purpose of this paper is to extend the optimal results on Berry–Esseen bounds for $U$-statistics of order 2 to those of higher orders. In stark contrast with studies on $U$-statistics of order 2, there is a very limited literature on optimal or near-optimal error bounds for $U$-statistics of higher orders. In this paper, we intend to fill in the gap. The work is of both theoretical and practical value since many symmetric statistics may be approximated arbitrarily closely by $U$-statistics of sufficiently high order under appropriate conditions. As an application, we will derive a near-optimal error bound for Studentized $U$-statistics of order 2 in Section 3.3.

At first sight, it appears to be an easy task to extend the optimal Berry–Esseen bounds from second-order $U$-statistics to higher order ones. (At least we naively thought so in the beginning.) However, a close inspection soon leads one to believe that this is far from trivial. We note that the usual approach of dealing with $U$-statistics is to first use the *Hoeffding decomposition* to turn the statistic of interest into a sum of sequentially smaller and uncorrelated terms, and then to use *truncation techniques* to each of these terms. Broadly speaking, there are two main difficulties with this approach which we must overcome. Specifically, we have the following goals in mind.

(A) *A direct and transparent methodology is needed.*
The truncation techniques are manageable for second-order $U$-statistics. However, as the order gets larger, it becomes more and more intangible and unworkable. *In order to treat $U$-statistics of higher order than 2, in this paper, we may have to abandon the idea of truncation, and instead choose a more direct and more transparent approach.*

(B) *A simple and unified form of error bounds is needed.*
For second-order $U$-statistics, the Hoeffding decomposition produces a linear statistic $\mathbb{L}$ plus a degenerate $U$-statistic $\mathbb{T}$ in (1.2). As a result, the two terms involving $E|g_1|^3$ and $E|\eta_{12}|^{5/3}$ in (1.3) are from the linear term $\mathbb{L}$ and the degenerate $U$-statistic term $\mathbb{T}$, respectively. More generally, applying the Hoeffding decomposition to $U$-statistics of order $k > 2$ leads to a sum of $k$ terms. Consequently, the resulting optimal error bound would contain $k$ terms which would become more complicated as $k$ gets larger. *In this paper, we strive to provide simple and general error bounds for $U$-statistics of general orders.*

Let us now look at (B) in more detail. Again, we will examine the case of $U$-statistics of order 2 first. A closer look at the proof of (1.3), as in [14] and [9], shows that a more refined Berry–Esseen bound can be given as follows:

$$(1.4) \qquad \|F_{\mathbb{L}+\mathbb{T}} - \Phi\| \ll \frac{1}{\sqrt{n}} \bigg( \frac{E|g_1|^3}{\sigma_g^3} + \frac{E|\eta_{12}|^{5/3}}{\sigma_g^{5/3}} + \frac{E|g_1 g_2 \eta_{12}|}{\sigma_g^3} \bigg).$$

Using some truncation arguments, it can be shown that (1.4) is equivalent to

$$(1.5) \qquad \|F_{\mathbb{L}+\mathbb{T}} - \Phi\| \ll \frac{1}{\sqrt{n}} \frac{E|g_1|^3}{\sigma_g^3} + \frac{1}{n} \frac{E|\eta_{12}|^2}{\sigma_g^2} + \frac{1}{\sqrt{n}} \frac{E|g_1 g_2 \eta_{12}|}{\sigma_g^3}.$$

Note that the last term in (1.5) shows the interaction effect between $\mathbb{L}^2$ and $\mathbb{T}$ in (1.2). This suggests that if we take the correlation between $\mathbb{L}$ and $\mathbb{T}$ into account, we might be able to improve the error bound. Indeed, if we let

$$N_2(x) = \Phi(x) + \frac{1}{\sqrt{n}} \kappa_2 \Phi'''(x) \qquad \text{where } \kappa_2 = E g_1 g_2 \eta_{12} / \sigma_g^3$$

be the *adjusted* normal approximation, Alberink and Bentkus [2, 3] show that the following optimal error bound holds:

$$\|F_{\mathbb{L}+\mathbb{T}} - N_2\| \ll \frac{1}{\sqrt{n}} \frac{E|g_1|^3}{\sigma_g^3} + \frac{1}{n} \frac{E|\eta_{12}|^2}{\sigma_g^2} \tag{1.6}$$

(see Theorem 1 of [2, 3]). We reiterate that the error bound given in (1.5) is optimal, and implies the optimal bound (1.3) as a special case by a simple truncation argument.

Finally, we can rewrite (1.5) in a very simple form

$$\|F_{\mathbb{L}+\mathbb{T}} - N_2\| \ll \|F_{\mathbb{L}} - \Phi\| + \operatorname{var} \mathbb{T}. \tag{1.7}$$

In other words, the error bounds in the adjusted normal approximations for $U$-statistics of order 2 are simply the error bounds for the dominant linear part $\mathbb{L}$ plus the variance of the remaining "error" term $\mathbb{T}$.

It is somewhat surprising that the optimal error bound for $U$-statistics of order 2 can be written in such a simple and compact form given by (1.7). The appearance of the variance term $\operatorname{var} \mathbb{T}$ in (1.7) is appealing since the variance is the most natural and easy-to-interpret measure to describe the effects of the error term $\mathbb{T}$. Furthermore, from a computational point of view, it is easier to calculate the variance of $\mathbb{T}$ than some other moments of $\mathbb{T}$.

The purpose of this paper is to derive optimal or near-optimal error bounds in the normal approximations for $U$-statistics of order $k \geq 2$, similar to those given in (1.7). We emphasize that the truncation techniques are heavily used in [2, 3] for second-order $U$-statistics, and will not work for $U$-statistics of higher orders.

The paper is arranged as follows. In Section 2, some definitions and notation will be given. In Section 3, we establish optimal or near-optimal error bounds for $U$-statistics of order $k \geq 2$. In Section 4, we give some approximations to expectations of smooth functions of $U$-statistics and some of their applications. Possible extensions to other statistics are discussed in Section 5. We provide the details of the proofs in Section 6.

**2. Definitions and notation.** Let $X, X_1, \ldots, X_n$ be a sequence of i.i.d. r.v.'s. Suppose that a symmetric statistic $\mathbb{S} = \mathbb{S}(X_1, \ldots, X_n)$ can be decomposed into

$$\mathbb{S} = \mathbb{L} + \mathbb{T}, \tag{2.1}$$

where $\mathbb{L}$ and $\mathbb{T}$ are symmetric $U$-statistics, $\mathbb{L}$ of order 1 and $\mathbb{T}$ of order $k$, respectively. This means that $\mathbb{L}$ is of the form

$$\mathbb{L} = \sum_{i=1}^{n} L_i, \tag{2.2}$$

and $\mathbb{T}$ is of the form, by way of Hoeffding decomposition,

(2.3)
$$\mathbb{T} = \sum_{i=1}^{n} T_i + \sum_{1 \leq i < j \leq n}^{n} T_{ij} + \cdots + \sum_{1 \leq i_1 < \cdots < i_k \leq n}^{n} T_{i_1 \cdots i_k}$$
$$\stackrel{\text{def}}{=} \mathbb{T}_1 + \cdots + \mathbb{T}_k,$$

where $L_i = l(X_i)$ and $T_{i_1 \cdots i_k} = t_k(X_{i_1}, \ldots, X_{i_k})$ are Borel measurable functions, and $t_k$ is invariant under permutation of variables. Without loss of generality, it is assumed that

(2.4) $$\mathbb{E}\mathbb{L} = \mathbb{E}\mathbb{T} = 0 \quad \text{and} \quad \text{var}\,\mathbb{L} = 1,$$

which can always be achieved after re-centering and re-scaling the relevant statistics. Typically, $\mathbb{L}$ is the dominant term and $\mathbb{T}$ can then be called the *error* term.

Now let us define, for $\alpha \in [1,2]$,

$$\beta = n\mathbb{E}|L_1|^3, \qquad \gamma = \text{var}\,\mathbb{T}, \qquad \gamma^{(\alpha)} = \sum_{p=1}^{k} \binom{n}{p} \mathbb{E}|T_{1 \cdots p}|^{\alpha}.$$

Note that

$$\mathbb{E}|\mathbb{L}|^3 \ll n\mathbb{E}|L_1|^3 = \beta,$$

(2.5)
$$\mathbb{E}|\mathbb{T}|^2 = \sum_{i=1}^{k} \mathbb{E}\mathbb{T}_i^2 = \sum_{i=1}^{k} \binom{n}{i} \mathbb{E}|T_{1 \cdots i}|^2 = \gamma,$$

$$\mathbb{E}|\mathbb{T}|^{\alpha} \leq k^{\alpha-1} \sum_{i=1}^{k} \mathbb{E}|\mathbb{T}_i|^{\alpha} \ll \sum_{i=1}^{k} \binom{n}{i} \mathbb{E}|T_{1 \cdots i}|^{\alpha} = \gamma^{(\alpha)},$$

where the second inequality in (2.5) follows by first applying an inequality (2.18) from [10],

$$\mathbb{E}|\mathbb{T}_i|^{\alpha} \leq C(\alpha, i)\mathbb{E}\left| \sum_{1 \leq j_1 < \cdots < j_i \leq n} T_{j_1 \cdots j_i}^2 \right|^{\alpha/2}$$

with $C(\alpha, i)$ a constant depending only on $\alpha$ and $i$, and then applying the following simple inequality (noting $\alpha/2 \leq 1$):

$$\left| \sum_{1 \leq j_1 < \cdots < j_i \leq n} T_{j_1 \cdots j_i}^2 \right|^{\alpha/2} \leq \sum_{1 \leq j_1 < \cdots < j_i \leq n} |T_{j_1 \cdots j_i}|^{\alpha}.$$

Therefore, $\beta$, $\gamma$ and $\gamma^{(\alpha)}$ are closely related to $E|\mathbb{L}|^3$, $\mathbb{E}\mathbb{T}^2$ and $\mathbb{E}|\mathbb{T}|^{\alpha}$, respectively. Furthermore, $\gamma = \gamma^{(2)}$.

Finally, the following convention will be adopted throughout the paper: conditions appearing in a statement are implied implicitly in that statement. For instance, in (3.4), the conditions $\beta < \infty$ and $\gamma < \infty$ are assumed even though (3.4) still holds true without these moment conditions.

**3. Error bounds for $U$-statistics.** Asymptotic normality of the linear term $\mathbb{L}$ and its rates of convergence are well known. For instance, the latter is described by the Berry–Esseen bound:

$$\|F_{\mathbb{L}} - \Phi\| \ll \beta. \tag{3.1}$$

If $\mathbb{T} \to 0$ in probability, we have $\|F_{\mathbb{S}} - \Phi\| \to 0$ as $n \to \infty$. Following the discussion in the Introduction, we are interested in the error bound for $\|F_{\mathbb{S}} - \Phi\|$ in the form of

$$\|F_{\mathbb{S}} - N_p\| \ll \|F_{\mathbb{L}} - \Phi\| + f(\mathbb{E}|\mathbb{T}|^\alpha) \ll \beta + f(\gamma^{(\alpha)}) \tag{3.2}$$

for some function $f$, $N_p$ for some $p \geq 0$ (to be defined below), and $\alpha \leq 2$. In particular, when $\alpha = 2$, (3.2) becomes

$$\|F_{\mathbb{S}} - N_p\| \ll \|F_{\mathbb{L}} - \Phi\| + f(\operatorname{var} \mathbb{T}) \ll \beta + f(\gamma).$$

In other words, we wish to investigate the question: *if a linear statistic $\mathbb{L}$ is perturbed by some error term $\mathbb{T}$, how do we correct the usual Berry–Berry bound for normal approximations, using only the variance of $\mathbb{T}$?*

The *adjusted normal approximation $N_p$* [which was mentioned in (3.2)] can be defined as follows. Define $N_0(x) = \Phi(x)$ and

$$N_p(x) = \Phi(x) + G_1(x) + \cdots + G_p(x), \qquad p = 1, 2, \ldots, \tag{3.3}$$

where, for $p = 1, \ldots, k$,

$$\kappa_p = \operatorname{cov}(\mathbb{L}^p, \mathbb{T}) = \mathbb{E}\mathbb{L}^p\mathbb{T} = \binom{n}{p}\mathbb{E}L_1 \cdots L_p\mathbb{T},$$

$$G_p(x) = (-1)^{p+1}\kappa_p\Phi^{(p+1)}(x).$$

Two questions immediately arise from the above definition:

1. *Under what conditions can we guarantee the existence of $\kappa_p$ for $p \geq 1$?*
   By the Hölder inequality, we have

   $$\mathbb{E}|L_1 \cdots L_i \mathbb{T}| \leq (\mathbb{E}|L_1|^3 \cdots \mathbb{E}|L_i|^3)^{1/3}(\mathbb{E}|\mathbb{T}|^{3/2})^{2/3}.$$

   Therefore, $\kappa_p$ exists if $\mathbb{E}|\mathbb{T}|^{3/2} < \infty$ as we assume that $\mathbb{E}|L_i|^3 < \infty$ for $1 \leq i \leq n$. If $\mathbb{E}|T|^{3/2} = \infty$, we will see later that we only use $N_0(x) = \Phi(x)$.
2. *How does one interpret $N_p$?*
   The adjusted normal approximation $N_p$ plays a central role in this paper. For $U$-statistics of order 2 which first appeared in [2], we have $\kappa_1 = 0$ and $\kappa_2 = -n^{-1/2}\mathbb{E}g_1g_2\eta_{12}/\sigma_g^3$, hence

   $$N_2(x) = \Phi(x) - \frac{1}{\sqrt{n}}\frac{\mathbb{E}g_1g_2\eta_{12}}{\sigma_g^3}\Phi'''(x).$$

Note that $N_2$ is different from its corresponding Edgeworth expansion,

$$F_{\mathbb{L}+\mathbb{T}}(x) = \Phi(x) - \frac{1}{\sqrt{n}}\left(\frac{\mathbb{E}g_1 g_2 \eta_{12}}{2\sigma_g^3} + \frac{\mathbb{E}g_1^3}{6\sigma_g^3}\right)\Phi'''(x) + o(n^{-1/2}),$$

which holds under the following optimal conditions: (i) the distribution of $g(X_1)$ is nonlattice; (ii) $\mathbb{E}|g_1|^3 < \infty$ and $\mathbb{E}|\eta_{12}|^{5/3} < \infty$ (see [12]). On the other hand, $N_p$'s and Edgeworth expansions are similar in the sense that both try to improve the normal approximations by including some higher order correction terms. We point out that Edgeworth expansions do require some smoothness conditions while Berry–Esseen bounds do not.

Of course, the adjusted normal approximations $N_p$ are not introduced to compete with Edgeworth expansions. Rather, they are convenient tool for obtaining a neat and unified error bound in the normal approximation for $U$-statistics of general order.

In the following, we will give error bounds under two different sets of conditions: $\mathbb{E}\mathbb{T}^2 < \infty$, and $\mathbb{E}|\mathbb{T}|^\alpha < \infty$ for $\alpha \in [1,2)$, respectively. Since the Hoeffding decomposition of $\mathbb{T}$ requires $\mathbb{E}|\mathbb{T}| < \infty$, we will not consider the case $\mathbb{E}|\mathbb{T}|^\alpha < \infty$ for $\alpha \in (0,1)$.

3.1. *Error bounds under $\mathbb{E}\mathbb{T}^2 < \infty$.* Here is our key result.

THEOREM 1. *We have*

(3.4) $$\|F_\mathbb{S} - N_k\| \ll \beta + \gamma \ln^2(\gamma),$$

*where the last term $\gamma \ln^2(\gamma)$ is taken to be zero when $\gamma = 0$.*

For second-order $U$-statistics, Alberink and Bentkus [2, 3] show

$$\|F_\mathbb{S} - N_k\| \ll \beta + \gamma,$$

which is sharper than that in Theorem 1 when $k = 2$. We conjecture that the logarithmic factor in Theorem 1 can be removed for $U$-statistics of order $k > 2$. Actually, we have derived an error bound without log factors which, for the $U$-statistics of order 2, implies the result of [2, 3]. Unfortunately, for $k \geq 3$, the error bound involves certain *conditional* variances, and the proof of this bound is indeed very complicated, compared to the proof of Theorem 1. We do not provide this bound here since, in applications, this bound has only minor advantages, compared to the simple and indeed convenient bound of Theorem 1.

Applying Theorem 1, we can derive the following optimal error bound for $\|F_\mathbb{S} - \Phi\|$ for $U$-statistics of order $k$. (A similar theorem to (3.5) is also derived in [7], using Stein's method.)

THEOREM 2. *We have*

(3.5) $$\|F_\mathbb{S} - \Phi\| \ll \beta + \gamma^{1/2}.$$

*Furthermore, the exponent $1/2$ in (3.5) is the best possible for $U$-statistics of order $k \geq 2$.*

Comparing Theorems 1 and 2, we see that the adjusted normal approximation $N_k(x)$, which takes into account of the correlations between $\mathbb{L}^k$ and $\mathbb{T}$, does show improvement over the standard normal distribution $\Phi(x)$.

3.2. *Error bounds under $\mathbb{E}|\mathbb{T}|^\alpha < \infty$ for $\alpha \in [1, 2)$.* Theorem 1 contains a log factor which will be inherited when applying the result. We will provide an optimal error bound without the log factor under a lower moment assumption $\mathbb{E}|\mathbb{T}|^\alpha < \infty$ for $\alpha \in [1, 2)$.

To do this, we define $N_0 = \Phi$ and

$$N_p = \Phi + G_1 + \cdots + G_p,$$

where $p = p(\alpha, k)$ is the largest integer such that $p < (\alpha - 1)/(2 - \alpha)$ and $p \leq k$. For example, if $\alpha \leq 3/2$, then $N_p = \Phi$. If $3/2 < \alpha \leq 5/3$ and $k \geq 1$, then $N_p = \Phi + G_1$. If $5/3 < \alpha \leq 7/4$ and $k \geq 2$, then $N_p = \Phi + G_1 + G_2$.

THEOREM 3. *For $\alpha \in [1, 2)$, let $p = p(\alpha, k)$ be defined above, we have*

$$\|F_\mathbb{S} - N_p\| \ll \beta + \gamma^{(\alpha)}.$$

Theorem 3 would be of most interest when $\alpha$ is close to 2 in applications. As an example, we will apply Theorem 3 below to obtain a near-optimal Berry–Esseen bound for Studentized $U$-statistics of order 2.

3.3. *An application to Studentized $U$-statistics of order 2.* The optimal error bound for standardized $U$-statistics of order 2 was given in (1.3). We conjecture that the same optimal error bound applies to their corresponding Studentized $U$-statistics defined by

$$\widehat{\mathbb{S}} = \frac{\sqrt{n}(U_n - \theta)}{2\widehat{\sigma}_g},$$

where $\widehat{\sigma}_g^2$ is the jackknife variance estimator,

$$\widehat{\sigma}_g^2 = \frac{n-1}{(n-2)^2} \sum_{i=1}^n \left( \frac{1}{n-1} \sum_{j=1, j \neq i}^n h_{ij} - U_n \right)^2.$$

However, we have so far only managed to prove a near-optimal result, as given in Theorem 4 below. Although this result is the best available bound in the literature, it falls a little short of the optimal bound (i.e., $\epsilon = 0$).

THEOREM 4. *For any $\epsilon > 0$, we have*

$$\|F_{\widehat{\mathbb{S}}} - \Phi\| \ll \frac{1}{\sqrt{n}}\bigg(\frac{E|g_1|^3}{\sigma_g^3} + \frac{E|\eta_{12}|^{5/3+\epsilon}}{\sigma_g^{5/3+\epsilon}}\bigg).$$

The proof of Theorem 4 is very involved, and hence will not be given here. Its proof can be obtained from the authors. The proof basically involves approximating $\widehat{\mathbb{S}}$ by a $U$-statistic of sufficiently high order and then applying nontrivial truncation techniques. Currently we are trying to eliminate $\epsilon$ in Theorem 4, and hope to be able to report on this in the near future.

**4. Error bounds for expectations of smooth functions.** Here, we give some approximations to expectations of smooth functions of $U$-statistics. These results (e.g., Theorems 5 and 6) are of interest in their own right. There are also other useful applications, for example, in approximating characteristic functions. For instance, Corollaries 1 and 2 have been used in the proofs of Theorems 1 and 3, respectively [see (6.6) and (6.22)].

Throughout this section, let $f:\mathbb{R} \to \mathbb{R}$ be a sufficiently smooth function.

THEOREM 5. *We have*

$$\bigg\|\int_\mathbb{R} f(u)\,d\mathbb{P}\{\mathbb{S} \le u\} - \int_\mathbb{R} f(u)\,dN_k(u)\bigg\| \ll C(\beta + \gamma),$$

*where $C \stackrel{\mathrm{def}}{=} \|f''\| + \cdots + \|f^{(k+4)}\|$.*

Choosing $f(u) = \exp\{itu\}$ in Theorem 5, we obtain a useful inequality for the characteristic functions.

COROLLARY 1. *We have*

$$\bigg|\mathbb{E}\exp\{it\mathbb{S}\} - \int_\mathbb{R}\exp\{itu\}\,dN_k(u)\bigg| \ll (|t| + |t|^{k+4})(\beta + \gamma).$$

Theorem 5 and Corollary 1 can be extended to the case of the lower moment assumption $\mathbb{E}|\mathbb{T}|^\alpha < \infty$ with $\alpha < 2$. However, the technical details are more involved in Theorem 6 than those in Theorem 5.

THEOREM 6. *For $\alpha \in [1,2)$, we have*

$$\bigg\|\int_\mathbb{R} f(u)\,d\mathbb{P}\{\mathbb{S} \le u\} - \int_\mathbb{R} f(u)\,dN_p(u)\bigg\| \ll C(\beta + \gamma^{(\alpha)}),$$

*where $p = p(\alpha, k)$ is defined in Theorem 3, and $C = \|f'\| + \|f'''\|$ for $\alpha = 1$ and $C = \|f'\| + \cdots + \|f^{(k+4)}\|$ for $\alpha \in (1,2)$.*

Choosing $f(u) = \exp\{itu\}$ in Theorem 6, we get:

COROLLARY 2. *For $\alpha \in [1,2)$, we have*

$$\left| \mathbb{E} \exp\{it\mathbb{S}\} - \int_{\mathbb{R}} \exp\{itu\} \, dN_p(u) \right| \ll C(t)(\beta + \gamma^{(\alpha)}),$$

*where $p = p(\alpha, k)$ is defined in Theorem 3, and $C(t) = |t| + |t|^3$ for $\alpha = 1$ and $C(t) = |t| + |t|^{k+4}$ for $\alpha \in (1,2)$.*

**5. Extensions.** One may consider extending the work on the error bounds from $U$-statistics to more general class of statistics. Consider a nonlinear statistic $\mathbb{W} = \mathbb{W}(X_1, \ldots, X_n)$, based on a sequence of the r.v.'s $X_1, \ldots, X_n$. Examples include symmetric statistics, $U$-statistics with non-i.i.d. observations, $L$-statistics, Studentized statistics, exchangeable statistics, finite population statistics, and so on. Error bounds for symmetric statistics and non-symmetric statistics were considered by [17] and [7]. Alberink [1] studied error bounds for $U$-statistics with non-i.i.d. observations.

Suppose that we can linearize $\mathbb{W}$ into

$$\mathbb{W} = \mathbb{L} + \Delta,$$

where $\mathbb{L}$ is a linear statistic and $\Delta$ is an error term. One way to linearize $\mathbb{W}$ is by the Hoeffding decomposition. We have seen how it works for $U$-statistics in this paper. For symmetric statistics, we refer to [17]. The Hoeffding decompositions for independent observations and nonsymmetric statistics was discussed by [13]. For orthogonal decomposition of finite population statistics, see [6]. Hoeffding decompositions for exchangeable statistics were considered by [16].

If $\mathbb{L}$ is asymptotically normal and $\Delta \to 0$ in probability, $\mathbb{W}$ is asymptotically normal. Using truncation and the Chebyshev inequality, we can easily get the following error bound:

$$\|F_{\mathbb{W}} - \Phi\| \leq \|F_{\mathbb{L}} - \Phi\| + (1+p)[(\|\Phi'\|)/p]^{p/(1+p)} (E|\Delta|^p)^{1/(1+p)}.$$

See (1.3) of [7]. By taking $p = 2$, we get

(5.1) $$\|F_{\mathbb{W}} - \Phi\| \leq \|F_{\mathbb{L}} - \Phi\| + 2(\operatorname{var} \Delta)^{1/3}.$$

For symmetric statistics, the next example shows that the exponent $1/3$ in (5.1) is already the best possible statistic. The example is similar to an example of [7].

EXAMPLE 1. *Take $X, X_1, \ldots, X_n$ to be i.i.d. $N(0,1)$ r.v.'s. Define*

(5.2) $\quad \mathbb{L} = n^{-1/2}(X_1 + \cdots + X_n), \qquad \Delta = -\varepsilon(|\mathbb{L}|^{-a} - \mathbb{E}|\mathbb{L}|^{-a}),$

where $0 < a < 1/2$ and $\varepsilon > 0$. Then $\mathbb{W} = \mathbb{L} + \Delta$ is a symmetric statistic with finite variance since $\mathbb{E}\Delta^2 \leq 2\varepsilon^2 \mathbb{E}|\mathbb{L}|^{-2a} < \infty$. Now if we have

$$\|F_\mathbb{W} - \Phi\| \leq c\|F_\mathbb{L} - \Phi\| + c(\operatorname{var}\Delta)^\vartheta \tag{5.3}$$

or

$$\|F_\mathbb{W} - \Phi\| \leq c\mathbb{E}|X|^3/\sqrt{n} + c(\operatorname{var}\Delta)^\vartheta, \tag{5.4}$$

then the exponent $\vartheta$ must satisfy $\vartheta \leq 1/3$.

PROOF. Clearly, $\mathbb{L} \sim N(0,1)$. So $\|F_\mathbb{L} - \Phi\| = \sup_x |\mathbb{P}\{\mathbb{L} \leq x\} - \Phi(x)| = 0$. Thus it suffices to show that the weaker inequality (5.4) implies $\vartheta \leq 1/3$. Also note that $\operatorname{var}\Delta = c_1\varepsilon^2$ where $c_1 = \operatorname{var}(|\mathbb{L}|^{-a})$ is some positive constant depending only on $a$. It follows that $\|F_\mathbb{W} - \Phi\|$ and $\operatorname{var}\Delta$ in this example do not depend on $n$. Letting $n \to \infty$ in (5.4) yields

$$\|F_\mathbb{W} - \Phi\| \leq c(\operatorname{var}\Delta)^\vartheta = cc_1^\vartheta \varepsilon^{2\vartheta}. \tag{5.5}$$

On the other hand, for sufficiently small $\varepsilon > 0$, we have

$$\begin{aligned}
\|F_\mathbb{W} - \Phi\| &= \sup_x |\mathbb{P}\{\mathbb{W} \leq x\} - \Phi(x)| \\
&\geq \mathbb{P}\{\mathbb{W} \leq \varepsilon\mathbb{E}|\mathbb{L}|^{-a}\} - \mathbb{P}\{\mathbb{L} \leq \varepsilon\mathbb{E}|\mathbb{L}|^{-a}\} \\
&= \mathbb{P}\{\mathbb{L}|\mathbb{L}|^a \leq \varepsilon\} - \mathbb{P}\{\mathbb{L} \leq \varepsilon\mathbb{E}|\mathbb{L}|^{-a}\} \\
&= \mathbb{P}\{0 \leq \mathbb{L}^{a+1} \leq \varepsilon\} - \mathbb{P}\{0 \leq \mathbb{L} \leq \varepsilon\mathbb{E}|\mathbb{L}|^{-a}\} \\
&\geq c_2\varepsilon^{1/(a+1)} - c_3\varepsilon,
\end{aligned} \tag{5.6}$$

where $c_2$ and $c_3$ are positive constants depending only on $a$. The inequalities (5.5) and (5.6) imply $cc_1^\vartheta \varepsilon^{2\vartheta} \geq c_2\varepsilon^{1/(a+1)} - c_3\varepsilon$, that is,

$$cc_1^\vartheta \varepsilon^{2\vartheta - 1/(a+1)} \geq c_2 - c_3\varepsilon^{a/(a+1)}. \tag{5.7}$$

If $2\vartheta - 1/(a+1) > 0$, then letting $\varepsilon \downarrow 0$ in (5.7) would imply that $0 \geq c_2 > 0$, a contradiction. Hence we must have $2\vartheta - 1/(a+1) \leq 0$ for all $0 < a < 1/2$. Letting $a \uparrow 1/2$, we have $2\vartheta - 2/3 \leq 0$, i.e., $\vartheta \leq 1/3$. □

**6. Proofs of Theorems 1–6.** We first prove several useful lemmas.

LEMMA 1. *We have:*

(a) $1 \leq \sqrt{n}\beta$;
(b) $n\mathbb{E}|L_1|^q \leq \beta^{q-2}$ for $2 \leq q \leq 3$;
(c) $|\kappa_s| \leq \sqrt{\gamma_s} \leq \sqrt{\gamma}$ for $s = 1, \ldots, k$;
(d) $|\kappa_s| \leq \beta^{\delta s} + \gamma_s^{(\alpha)} \leq \beta^{\delta s} + \gamma^{(\alpha)}$ where $3/2 < \alpha < 2$, $\delta = (2-\alpha)/(\alpha-1)$, and $s$ is an integer such that $1 \leq s < \delta^{-1} \wedge k$.

PROOF. We will use Hölder's inequality: $\mathbb{E}|XY| \le |(\mathbb{E}|X|^a)^{1/a}(\mathbb{E}|Y|^b)^{1/b}$ with $a > 1$ and $1/a + 1/b = 1$. Recall that $\mathbb{E}L^2 = n\mathbb{E}L_1^2 = 1$.

(a) $1 = (nn^{-1})^{3/2} = (n\mathbb{E}L_1^2)^{3/2} = n^{3/2}(\mathbb{E}L_1^2)^{3/2} \le n^{3/2}\mathbb{E}|L_1|^3 = \sqrt{n}\beta$.
(b) $n\mathbb{E}|L_1|^q = n\mathbb{E}|L_1|^{2(3-q)}|L_1|^{3(q-2)} \le n(\mathbb{E}L_1^2)^{3-q}(\mathbb{E}|L_1|^3)^{q-2} = \beta^{q-2}$.
(c) $|\kappa_s| \le \binom{n}{s}\mathbb{E}|L_1 \cdots L_s T_{1\cdots s}| \le \binom{n}{s}(\mathbb{E}L_1^2 \cdots \mathbb{E}L_s^2 \mathbb{E}T_{1\cdots s}^2)^{1/2} = (\binom{n}{s}/n^s)^{1/2}\gamma_s^{1/2} \le \gamma_s^{1/2} \le \gamma^{1/2}$.
(d) We use Hölder's inequality with exponents $a = \alpha$ and $b = \alpha/(\alpha - 1)$ and another basic inequality $x^{1/a}y^{1/b} \le x + y$ for $x, y \ge 0$. Also note that $n\mathbb{E}|L_1|^q \le \beta^{q-2}$ from (b), and $q - 2 = \delta$. We have

$$|\kappa_s| \le \binom{n}{s}\mathbb{E}|L_1 \cdots L_s T_{1\cdots s}| \le \binom{n}{s}(\mathbb{E}|L_1|^q)^{s/q}(\mathbb{E}|T_{1\cdots s}|^\alpha)^{1/\alpha}$$

$$\le \binom{n}{s}(\mathbb{E}|L_1|^q)^s + \binom{n}{s}\mathbb{E}|T_{1\cdots s}|^\alpha \le (n\mathbb{E}|L_1|^q)^s + \gamma_s^{(\alpha)}$$

$$\le \beta^{(q-2)s} + \gamma_s^{(\alpha)} = \beta^{\delta s} + \gamma_s^{(\alpha)}. \qquad \square$$

LEMMA 2. *For sufficiently smooth functions $g: \mathbb{R} \to \mathbb{R}$, one has*

(6.1) $$\|g(x+h) - g(x)\| \le (\|g\| + \|g'\|)|h|^{1/2}.$$

PROOF. The result follows from multiplying the obvious inequalities

$$\|g(x+h) - g(x)\| \le 2\|g\| \quad \text{and} \quad \|g(x+h) - g(x)\| \le \|g'\| |h|,$$

and taking a square root, and then applying $\sqrt{2ab} \le a + b$. $\square$

6.1. *Proof of Theorem 1.* Denote $\Delta \stackrel{\text{def}}{=} \|F_\mathbb{S} - N_k\|$. We consider two separate cases:

Case I: $\max\{\beta, \gamma_1, \ldots, \gamma_k\} > c_k$;

Case II: $\max\{\beta, \gamma_1, \ldots, \gamma_k\} \le c_k$, where $c_k > 0$ is a constant depending only on $k$ to be determined later.

Case I is relatively easy to prove. For case II, we take a classical approach. That is, using the well-known Esseen's smoothing inequality, we reduce the problem to the estimation of some characteristic functions.

PROOF FOR CASE I. We shall prove that $\Delta \ll \beta + \gamma$, which clearly implies the desired $\Delta \ll \beta + \gamma \ln^2(\gamma)$. From Lemma 1, we have $\kappa_p^2 \le \gamma_p$. Using this inequality, estimating $\mathbb{P}\{\mathbb{S} \le x\} \le 1$ and $|\Phi^{(p)}(x)| \ll 1$, we derive

(6.2) $\Delta \ll 1 + |\kappa_1| + \cdots + |\kappa_k| \ll 1 + \sqrt{\gamma_1} + \cdots + \sqrt{\gamma_k} \ll 1 + \sqrt{\gamma}$.

We consider the alternative cases $\gamma \ge 1$ and $\gamma < 1$ separately. If $\gamma \ge 1$, then using (6.2) we have $\Delta \ll \sqrt{\gamma} \le \gamma \le \beta + \gamma$. If $\gamma < 1$, then (6.2) implies $\Delta \ll 1$. Using the condition, $\max\{\beta, \gamma_1, \ldots, \gamma_k\} > c_k$, we have

$$\Delta \ll 1 \le c_k^{-1} \max\{\beta, \gamma_1, \ldots, \gamma_k\} \le c_k^{-1}(\beta + \gamma_1 + \cdots + \gamma_k) \ll \beta + \gamma.$$

PROOF FOR CASE II. We can assume that $n$ is sufficiently large so that

$$n \geq 1/c_k.$$

Indeed, for $n \leq 1/c_k$ the bound of Theorem 1 holds trivially. To see this, we note that, by (c) of Lemma 1, we have $\kappa_p^2 \leq \gamma_p$. Since $\max\{\beta, \gamma_1, \ldots, \gamma_k\} \leq c_k$ by assumption, we then have $\kappa_p^2 \leq \gamma_p \leq c_k \leq n^{-1} \leq 1$. Hence, it follows that

$$|N_k(x)| \ll 1 + |\kappa_1| + \cdots + |\kappa_k| \leq 1 + k \ll 1.$$

Therefore $\Delta \ll 1$. Using $1 \leq \sqrt{n}\beta$ (see Lemma 1), for small $n \leq 1/c_k$, we have $1 \ll \beta$. Inequalities $\Delta \ll 1$ and $1 \ll \beta$ imply $\Delta \ll \beta \ll \beta + \gamma$ which is somewhat better than the bound of Theorem 1.

By Lemma 1, $\kappa_p^2 \leq \gamma_p$ which combined with the condition, $\max\{\beta, \gamma_1, \ldots, \gamma_k\} \leq c_k$, implies that $\kappa_p^2 \ll 1$ for $1 \leq p \leq k$. Thus the function $N_k(x)$ has a bounded derivative,

(6.3) $$|N_k'(x)| \ll 1 + |\kappa_1| + \cdots + |\kappa_k| \ll 1.$$

Due to (6.3), to estimate $\Delta$ we can apply Esseen's smoothing inequality (see, e.g., [8], Chapter XVI, Lemma 3.2). For any $a > 0$ we have

(6.4) $$\Delta \ll \frac{1}{a} + \int_{-a}^{a} |t|^{-1} |f(t) - g(t)| \, dt,$$

where

(6.5) $$f(t) = \mathbb{E}\exp\{it\mathbb{S}\}, \qquad g(t) = \int_{\mathbb{R}} \exp\{itx\} \, dN_k(x).$$

Note that $g(t) = (1 + \kappa_1 (it)^3 + \cdots + \kappa_k (it)^{k+1}) \exp\{-t^2/2\}$. We choose

$$a = \frac{\sqrt{c_k}}{\beta + \gamma}.$$

Since $\max\{\beta, \gamma_1, \ldots, \gamma_k\} \leq c_k$, we have $\beta \leq c_k$ and $\gamma \leq kc_k$, resulting in

$$a \geq \sqrt{c_k}/(c_k + \cdots + c_k) = 1/((1+k)\sqrt{c_k}) > 1,$$

if $c_k$ is chosen small enough, for example, $c_k \leq (1+k)^{-2}$.

Split the integral in (6.4) as $\int_{-a}^{a} = \int_{|t|<C} + \int_{C<|t|<a}$, where $C = C_k$ is a sufficiently large positive constant depending only on $k$ to be chosen later.

To estimate $\int_{|t|<C}$ we use Corollary 1 of Section 4,

(6.6) $$|f(t) - g(t)| \ll (|t| + |t|^{k+4})(\beta + \gamma).$$

It follows that $\int_{|t|<C} \ll \beta + \gamma$, a bound which is somewhat better than the desired bound with $\gamma \ln^2(3 + 1/\gamma)$ in place of $\gamma$.

It remains to consider the integral $\int_{C \leq |t| \leq a}$. Introduce the characteristic functions $\vartheta = \vartheta(t)$ and $\varrho = \varrho(t)$,

(6.7) $$\vartheta = \mathbb{E}\exp\{itL_1\}, \qquad \varrho = \exp\{-t^2/(2n)\}.$$

By Lemma 3 below, we have, for $4k \leq m \leq n/4$,

(6.8) $\quad |f(t) - g(t)| \ll |t|^{k+4}(|\vartheta|^m + \varrho^m)(\beta + \gamma) + \gamma t^2 m/n \quad$ for $|t| \geq 1$.

We shall use the well-known simple bound

(6.9) $\quad\quad\quad\quad\quad |\vartheta| \leq \exp\{-t^2/(4n)\} \quad$ for $|t| \leq 1/\beta$.

In particular, the inequality (6.9) holds for $C < |t| < a$.

We choose the integer number $4k \leq m \leq n$ as

$$m = [4(k+5)nt^{-2}\ln|t|], \quad C < |t| < a,$$

where $[x]$ is the integer part of $x \in \mathbb{R}$. The number $m = m(t)$ depends on $t$. If $C = C_k$ is sufficiently large, then $m$ is a well defined integer such that $4k \leq m \leq n/4$, for sufficiently large $n$. Now (6.8) and (6.9) imply

(6.10) $\quad |f(t) - g(t)| \ll |t|^{-1}(\beta + \gamma) + \gamma \ln|t|, \quad C < |t| < a$.

Integrating (6.10) over $C < |t| < a$ we derive $\int_{C \leq |t| \leq a} \ll \beta + \gamma + \gamma \ln^2 a$. To conclude the proof of the theorem, we note that $1 \leq a \leq 3 + 1/\gamma$ due to our choices of constants, and therefore $\gamma + \gamma \ln^2 a \ll \gamma \ln^2(\gamma)$.

We now prove the following lemma, which is used in the proof above. The lemma gives an expansion of the characteristic functions for $|t| \geq 1$.

LEMMA 3. *Assume that $\beta \leq 1$ and $\gamma \leq 1$. Then for $n \geq 4k$ and $4k \leq m \leq n/4$ we have*

$$|f(t) - g(t)| \ll |t|^{k+4}(|\vartheta|^m + \varrho^m)(\beta + \gamma) + \gamma t^2 m/n \quad \text{for } |t| \geq 1.$$

PROOF. The proof follows from Propositions 1–3 below. When applying Propositions 1 and 2 one has to replace $m$ by $2m$. □

Let us introduce some notation first. Let $\Omega_m = \{1, \ldots, m\}$ and let $A$ be a subset of $\Omega_n = \{1, \ldots, n\}$, and use $|A|$ to denote the number of elements in $A$. For convenience, we write $T_{i_1,\ldots,i_k}$ or $T_{i_1\cdots i_k}$ instead of $T_{\{i_1,\ldots,i_k\}}$. Then we can write $\mathbb{T} = \mathbb{T}_1 + \cdots + \mathbb{T}_k$ where, for $1 \leq p \leq k$,

$$\mathbb{T}_p = \sum_{|A|=p} T_A \quad \text{with } T_A = t_p(X_j, j \in A).$$

We now split $\mathbb{T}$ into two parts $\mathbb{T}^{(m)}, \mathbb{T}^{(0)}$ so that

$$\mathbb{T} = \mathbb{T}^{(m)} + \mathbb{T}^{(0)}, \quad \mathbb{T}^{(m)} = \mathbb{T}_1^{(m)} + \cdots + \mathbb{T}_k^{(m)}, \quad \mathbb{T}^{(0)} = \mathbb{T}_1^{(0)} + \cdots + \mathbb{T}_k^{(0)},$$

where

$$\mathbb{T}_s^{(m)} = \sum_{|A|=s, A \cap \Omega_m \neq \varnothing} T_A \quad \text{and} \quad \mathbb{T}_s^{(0)} = \sum_{|A|=s, A \cap \Omega_m = \varnothing} T_A.$$

We are now ready to prove Propositions 1–3.

PROPOSITION 1. *Let $1 \leq m \leq n$. The characteristic function*
$$f = f(t) = \mathbb{E}\exp\{it\mathbb{S}\} \equiv \mathbb{E}\exp\{it(\mathbb{L} + \mathbb{T}^{(0)} + \mathbb{T}^{(m)})\}$$
*satisfies $f = f_0 + f_1 + R$ with*

(6.11) $\quad f_0 = \mathbb{E}\exp\{it(\mathbb{L} + \mathbb{T}^{(0)})\}, \qquad f_1 = (it)\mathbb{E}(\exp\{it(\mathbb{L} + \mathbb{T}^{(0)})\}\mathbb{T}^{(m)})$

*and a remainder term $R$ such that $|R| \ll \gamma t^2 m/n$.*

PROOF. We use $\tau, \tau_0, \tau_1, \tau_2, \ldots$ to denote i.i.d. r.v.'s which are uniformly distributed on $[0,1]$, and further assume that they are independent of all other r.v.'s. Given a smooth function $g$, we will frequently use the Taylor expansion of the following form
$$g(x+h) = \sum_{s=0}^{k} g^{(s)}(x)h^s/s! + \mathbb{E}(1-\tau)^k g^{(k+1)}(x+\tau h) h^{k+1}/k!.$$

For instance, we can expand $f$ in powers of $\mathbb{T}^{(m)}$, and obtain $f = f_0 + f_1 + R$ with
$$R = (it)^2 \mathbb{E}(1-\tau)\exp\{it(\mathbb{L} + \mathbb{T}^{(0)} + \tau\mathbb{T}^{(m)})\}(\mathbb{T}^{(m)})^2.$$
Thus
$$|R| \leq t^2 \mathbb{E}(\mathbb{T}^{(m)})^2 \ll t^2 \sum_{p=1}^{k} mn^{p-1}\mathbb{E}T_{1\cdots p}^2 = t^2 \frac{m}{n}\sum_{p=1}^{k} n^p \mathbb{E}T_{1\cdots p}^2 \ll \frac{m}{n}t^2\gamma. \quad \square$$

PROPOSITION 2. *Let $2k \leq m \leq n$. The function $f_0 + f_1$ with $f_0$ and $f_1$ defined by (6.11) satisfies $f_0 + f_1 = f_2 + f_3 + R_1$ with $f_2 = \mathbb{E}\exp\{it\mathbb{L}\}$, $f_3 = \mathbb{E}(it)\exp\{it\mathbb{L}\}\mathbb{T}$, and the remainder term $R_1$ satisfies $|R_1| \ll \gamma t^2 |\vartheta|^{m/2}$.*

PROOF. Since $\mathbb{T} = \mathbb{T}^{(0)} + \mathbb{T}^{(m)}$, it suffices to check that

(6.12) $\quad f_0 \sim f_2 + (it)\mathbb{E}\exp\{it\mathbb{L}\}\mathbb{T}^{(0)}, \qquad f_1 \sim (it)\mathbb{E}\exp\{it\mathbb{L}\}\mathbb{T}^{(m)}$

with remainder terms bounded as $R_1$.

Let us prove the first relation in (6.12). Note that
$$f_0 = \vartheta^m \mathbb{E}\exp\{it(\mathbb{L}^{(0)} + \mathbb{T}^{(0)})\} \qquad \text{where } \mathbb{L}^{(0)} = L_{m+1} + \cdots + L_n$$
is the part of $\mathbb{L}$ independent of $X_1, \ldots, X_m$. Now we can expand in powers of $\mathbb{T}^{(0)}$. We estimate the remainder term similar to the proof of Proposition 1. To estimate the variance of $\mathbb{T}^{(0)}$ we use the obvious inequality $\operatorname{var}\mathbb{T}^{(0)} \leq \operatorname{var}\mathbb{T}$.

Let us prove the second relation in (6.12). We consider only the case where $m$ is an even integer. It suffices to check that

(6.13) $\quad \mathbb{E}\exp\{it(\mathbb{L} + \mathbb{T}^{(0)})\}\mathbb{T}_p^{(m)} \sim \mathbb{E}\exp\{it\mathbb{L}\}\mathbb{T}_p^{(m)} \qquad \text{for } p = 1, \ldots, k.$

Let us show that

(6.14) $\quad \mathbb{E}\exp\{it(\mathbb{L}+\mathbb{T}^{(0)})\}\mathbb{T}_p^{(m)} = \mathbb{E}\exp\{it(\mathbb{L}+\mathbb{T}^{(0)})\}\mathbb{T}_p^*$

with $\mathbb{T}_p^* = \sum_{s=1}^p \sum_{(s)} \alpha_s T_A$ and $\alpha_s = \binom{m}{s}/\binom{m/2}{s}$, where $\sum_{(s)}$ is taken over all subsets $A \subset \Omega_n \setminus \Omega_{m/2}$ such that $|A \cap (\Omega_m \setminus \Omega_{m/2})| = s$.

To prove (6.14), let us start with a representation of $\mathbb{E}\exp\{it(\mathbb{L}+\mathbb{T}^{(0)})\}\mathbb{T}_p^{(m)}$. Sorting subsets $A$ according to the cardinality of the intersection $A \cap \{1,\ldots,m\}$ and using the symmetry and i.i.d. assumptions, we get

(6.15)
$$\begin{aligned}
&\mathbb{E}\exp\{it(\mathbb{L}+\mathbb{T}^{(0)})\}\mathbb{T}_p^{(m)} \\
&\quad = \sum_{A:\,|A|=p, A\cap\Omega_m \neq \varnothing} \mathbb{E} T_A \exp\{it(\mathbb{L}+\mathbb{T}^{(0)})\} \\
&\quad = \sum_{s=1}^p \sum_{A:\,|A|=p,|A\cap\Omega_m|=s} \mathbb{E} T_A \exp\{it(\mathbb{L}+\mathbb{T}^{(0)})\} \\
&\quad = \sum_{s=1}^p \binom{m}{s} \sum_{B:\,|B|=p-s,|B\cap\Omega_m|=\varnothing} \mathbb{E} T_{A_s \cup B} \exp\{it(\mathbb{L}+\mathbb{T}^{(0)})\},
\end{aligned}$$

where $A_1,\ldots,A_s \subset \{1,\ldots,m\}$ are arbitrary fixed subsets of cardinality $s$. For example, $A_s = \{1,\ldots,s\}$. A consideration similar to (6.15), starting with $\mathbb{E}\exp\{it(\mathbb{L}+\mathbb{T}^{(0)})\}\mathbb{T}_p^*$ instead of $\mathbb{E}\exp\{it(\mathbb{L}+\mathbb{T}^{(0)})\}\mathbb{T}_p^{(m)}$, shows that $\mathbb{E}\exp\{it(\mathbb{L}+\mathbb{T}^{(0)})\}\mathbb{T}_p^*$ is equal to the right-hand side of (6.15), which proves the identity (6.14).

The statistic $\mathbb{T}_p^*$ is independent of $X_1,\ldots,X_{m/2}$. Therefore, (6.14) implies

(6.16) $\quad \mathbb{E}\exp\{it(\mathbb{L}+\mathbb{T}^{(0)})\}\mathbb{T}_p^{(m)} = \vartheta^{m/2}\mathbb{E}\exp\{it(\mathbb{L}^{(m/2)}+\mathbb{T}^{(0)})\}\mathbb{T}_p^*.$

Now we can expand in powers of $\mathbb{T}^{(0)}$. This leads to

$$\mathbb{E}\exp\{it(\mathbb{L}+\mathbb{T}^{(0)})\}\mathbb{T}_p^{(m)} \sim \vartheta^{m/2}\mathbb{E}\exp\{it\mathbb{L}^{(m/2)}\}\mathbb{T}_p^* \equiv \mathbb{E}\exp\{it\mathbb{L}\}\mathbb{T}^{(m)}$$

up to an error bounded by $t^2|\vartheta|^{m/2}\mathbb{E}|\mathbb{T}^{(0)}\mathbb{T}_p^*|$ [this can be checked as we did for (6.15)]. Using the Hölder inequality, we have

$$\mathbb{E}|\mathbb{T}^{(0)}\mathbb{T}_p^*| \leq (\text{var}\,\mathbb{T}_p \,\text{var}\,\mathbb{T}_p^*)^{1/2} \ll \text{var}\,\mathbb{T} = \gamma,$$

since $\text{var}\,\mathbb{T}_p^* \ll \text{var}\,\mathbb{T}_p \leq \text{var}\,\mathbb{T}$. Combining the bounds, we get (6.13). □

PROPOSITION 3. *Assume that $|t| \geq 1$. Let $\beta \leq 1$ and $\gamma \leq 1$. Let $n \geq 4k$ and $1 \leq m \leq n/4$. Then the function $f_4 = \mathbb{E}\exp\{it\mathbb{L}\} + (it)\mathbb{E}\exp\{it\mathbb{L}\}\mathbb{T}$ satisfies $f_4 = g + R_2$ with $g(t) = (1 + \kappa_1(it)^3 + \cdots + \kappa_{k-1}(it)^{k+1})\varrho^n$ and the remainder term $R_2$ such that $|R_2| \ll (\beta+\gamma)|t|^{k+4}(\varrho^m + |\vartheta|^m).$*

PROOF. Using the i.i.d. and symmetry assumptions, we can write $f_4 = \vartheta^n + it \sum_{p=1}^{k} \vartheta^{n-p} \binom{n}{p} D_p$, where

$$D_p = \mathbb{E} T_{1 \cdots p} \exp\{it L_1\} \cdots \exp\{it L_p\}$$
$$= (it)^p \mathbb{E} L_1 \cdots L_p T_{1 \cdots p} \exp\{it(\tau_1 L_1 + \cdots + \tau_p L_p)\}.$$

In view of the inequalities $|\exp\{itu\} - 1| \ll \sqrt{|u|}$ and $\sqrt{u_1 + \cdots + u_p} \ll \sqrt{u_1} + \cdots + \sqrt{u_p}$ for $u_s \geq 0$, and $\sqrt{|t|} \leq |t|$ for $|t| \geq 1$, and the assumptions $n\mathbb{E}|L_1|^3 = \beta$ and $\mathbb{E} L_1^2 = 1/n$, we have

$$\left| \binom{n}{p} D_p - (it)^p \kappa_p \right| \leq |t|^p \binom{n}{p} \mathbb{E} |L_1 \cdots L_p T_{1 \cdots p}| \left| \exp\left\{it \left(\sum_{j=1}^{p} \tau_j L_j\right)\right\} - 1 \right|$$

$$\ll \binom{n}{p} |t|^{p+1/2} \mathbb{E} |L_1 \cdots L_p T_{1 \cdots p}| \left(\sum_{j=1}^{p} \sqrt{|L_j|}\right)$$

$$= p \binom{n}{p} |t|^{p+1/2} \mathbb{E} |L_1|^{3/2} |L_2 \cdots L_p T_{1 \cdots p}|$$

$$\ll |t|^{k+1} \sqrt{\beta} n^{-p/2} \binom{n}{p} (\mathbb{E} T_{1 \cdots p}^2)^{1/2}$$

$$\ll |t|^{k+1} \sqrt{\beta} \sqrt{\gamma}$$

$$\leq |t|^{k+1} (\beta + \gamma).$$

This in turn leads to $f_4 = f_5 + R_3$ with $f_5 = \vartheta^n + \sum_{p=1}^{k} (it)^{p+1} \kappa_p \vartheta^{n-p}$ and the remainder term $R_3$ such that $|R_3| \ll |t|^{k+2} |\vartheta|^m (\beta + \gamma)$.

Recall $1 \leq \sqrt{n}\beta$. By a simple Taylor expansion, we have $\beta \leq 1$, and $|\vartheta^p - 1| \ll t^2 \mathbb{E} L_1^2 = t^2/n \leq t^2 \beta^2 \leq t^2 \beta$. Hence $|\vartheta^{n-p} - \vartheta^n| \ll t^2 \beta |\vartheta|^m$. Therefore, using $\kappa_p \leq \sqrt{\gamma} \leq 1$, we can write $f_5 = f_6 + R_4$ with

$$f_6 = \vartheta^n \left(1 + \sum_{p=1}^{k} (it)^{p+1} \kappa_p\right)$$

and a remainder term $R_4$ such that $|R_4| \ll |t|^{k+3} |\vartheta|^m (\beta + \gamma)$. A very standard calculation shows that

(6.17) $$|\vartheta^n - \varrho^n| \ll |t|^3 (|\vartheta|^m + \varrho^m) \beta.$$

Thus using again $\kappa_p \leq \sqrt{\gamma} \leq 1$, it follows that $f_6 = g(t) + R_5$ with $|R_5| \ll |t|^{k+4} (|\vartheta|^m + \varrho^m) \beta$. Collecting all the inequalities, we complete the proof. $\square$

6.2. *Proof of Theorem 3.* The proof is similar to that of Theorem 1. Recall that $\gamma_s^{(\alpha)} = \binom{n}{s}\mathbb{E}|T_{1\cdots s}|^\alpha$. Without loss of generality we assume that

(6.18) $$\beta \leq c_k \quad \text{and} \quad \gamma_s^{(\alpha)} \leq c_k \quad \text{for } s = 1, \ldots, k;$$

and

(6.19) $$n \geq 1/c_k,$$

where $c_k$ is a sufficiently small positive constant. A reduction leading to (6.18) and (6.19) is based on applications of $|\kappa_s| \leq \beta^{\delta s} + \gamma_s^{(\alpha)} \leq \beta^{\delta s} + \gamma^{(\alpha)}$ (see Lemma 1) instead of $\kappa_p^2 \leq \gamma_p$ used in the proof of Theorem 1. Using (6.18), (6.19) and (d) of Lemma 1, it is easy to check that the function $N_p$ has a bounded derivative. Thus we can apply Esseen's smoothing inequality. For any $a > 0$,

(6.20) $$\|F_\mathbb{S} - N_p\| \ll \frac{1}{a} + \int_{-a}^{a} |t|^{-1}|f(t) - g^{(\alpha)}(t)|\,dt$$

with $f$ defined by (6.5) and

(6.21) $$g^{(\alpha)}(t) = \int_\mathbb{R} \exp\{itx\}\,dN_p(x).$$

We choose $a = \sqrt{c_k}/(\beta + \gamma^{(\alpha)})$. Split the integral in (6.20) as $\int_{-a}^{a} = \int_{|t|<C} + \int_{C<|t|<a}$, where $C = C_k$ is a sufficiently large positive constant depending only on $k$ and $\alpha$ to be chosen later.

To estimate $\int_{|t|<C}$ we use Corollary 2 of Section 4. Estimating $|t| \leq C$, the corollary implies

(6.22) $$|f(t) - g^{(\alpha)}(t)| \ll |t|^\alpha(\beta + \gamma^{(\alpha)}).$$

It follows that $\int_{|t|<C} \ll \beta + \gamma^{(\alpha)}$. Note that the presence of the factor $|t|^\alpha$ guaranties the convergence of the integral in a neighborhood of $t = 0$.

It remains to consider the integral $\int_{C \leq |t| \leq a}$. By Lemma 4 below, we have

$$|f(t) - g^{(\alpha)}(t)| \ll |t|^{k+4}(|\vartheta|^m + \varrho^m)(\beta + \gamma^{(\alpha)}) + \gamma^{(\alpha)}|t|^\alpha m/n$$

for $|t| \geq 1$ and $4k \leq m \leq n/4$. We integrate this bound as in the proof of Theorem 1. A small difference arises since now instead of $\gamma|t|^2 m/n$ we have the summand $\gamma^{(\alpha)}|t|^\alpha m/n$. The choice of $m \sim nt^{-2}\ln|t|$ leads to $\gamma^{(\alpha)}|t|^\alpha m/n \sim \gamma^{(\alpha)}|t|^{\alpha-2}\ln|t|$ which is an integrable function with respect to the measure $dt/|t|$ at $|t| = \infty$ since we assume that $\alpha < 2$. As a consequence, we now no longer have any log factors.

We now prove the following lemma, which is used in the proof above. This lemma extends Lemma 3 to the case of lower moment assumption.

LEMMA 4. *Let $1 \leq \alpha \leq 2$. Assume $\beta \leq 1$ and $\gamma \leq 1$. Then for $n \geq 4k$ and $4k \leq m \leq n/4$, we have*

$$|f(t) - g^{(\alpha)}(t)| \ll |t|^{k+4}(|\vartheta|^m + \varrho^m)(\beta + \gamma^{(\alpha)}) + \gamma^{(\alpha)}|t|^\alpha m/n \qquad \text{for } |t| \geq 1,$$

*where $\vartheta$, $\varrho$, $f$ and $g^{(\alpha)}$ are given by (6.7), (6.5) and (6.21), respectively.*

PROOF. We consider the cases $\alpha = 1$ and $1 < \alpha \leq 2$ separately. First consider the case $\alpha = 1$. Using $\exp\{itx\} = 1 + O(|tx|^\alpha)$, the expansion of Proposition 1 is replaced by $f = f_0 + R$ with a remainder $R$ such that $|R| \ll \gamma^{(\alpha)}|t|^\alpha m/n$. The function $f_0$ can be represented as

$$f_0 \equiv \mathbb{E}\exp\{it(\mathbb{L} + \mathbb{T}^{(0)})\} = \vartheta^m \mathbb{E}\exp\{it(\mathbb{L}^{(0)} + \mathbb{T}^{(0)})\}.$$

Thus by an expansion in powers of $\mathbb{T}^{(0)}$ we can replace it by $\mathbb{E}\exp\{it\mathbb{L}\} = \vartheta^n$. We conclude the proof replacing $\vartheta^n$ by $\varrho^n \equiv g \equiv g^{(\alpha)}$ [see (6.17)].

Next consider the case $1 < \alpha \leq 2$. Using $\exp\{itx\} = 1 + itx + O(|tx|^\alpha)$ and repeating arguments used in the proof of Proposition 1, we have $f = f_0 + f_1 + R$ with a remainder term $R$ such that $|R| \ll \gamma^{(\alpha)}|t|^\alpha m/n$. The functions $f_0, f_1$ are the same as in Proposition 1. A repetition of the proof of Proposition 2 leads to $f_0 + f_1 = f_2 + f_3 + R_1$ with $|R_1| \ll t^2 \gamma^{(\alpha)} |\vartheta|^{m/2}$. A small difference in the proof is that now we have to use Taylor expansions and Hölder's inequalities adjusted to our lower moment assumption.

To conclude the proof we have to show that $f_2 + f_3 = g^{(\alpha)} + R_2$ with a remainder term $R_2$ such that $|R_2| \ll (\beta + \gamma^{(\alpha)})t^{k+4}(\varrho^m + |\vartheta|^m)$. The proof is similar to that of Proposition 3 provided that some adjustments related to the lower moment assumption are made. The flavor of the adjustments is like that used in the proof of Theorem 6. We omit formal exposition of indeed lengthy and technical details. □

6.3. *Proof of Theorem 2.* We prove the first half now. If $\gamma \geq 1$, then we estimate $\|F_\mathbb{S} - \Phi\| \leq 2$, and it follows that $\|F_\mathbb{S} - \Phi\| \ll \beta + \gamma$. In the case of $\gamma \leq 1$, we use the bound of Theorem 1. Since $|G_s| \ll |\kappa_s|$ and $\kappa_s^2 \leq \gamma$ from Lemma 1, we have

$$\|F_\mathbb{S} - \Phi\| \leq \|F_\mathbb{S} - N_k\| + \sum_{s=1}^{k} |G_s| \ll \beta + \gamma \ln^2(\gamma) + \sqrt{\gamma} \ll \beta + \sqrt{\gamma}.$$

Next, we will show the second half of the theorem with an example. Let $X, X_1, \ldots, X_n$ be i.i.d. $N(0,1)$ r.v.'s, and define $\mathbb{S} = \mathbb{L} + \mathbb{T}$, where

$$\mathbb{L} = \frac{1}{\sqrt{n}} \sum_{i=1}^{n} X_i, \qquad \mathbb{T} = \frac{2\varepsilon}{(n-1)\sqrt{n}} \sum_{1 \leq i < j \leq n} X_i X_j$$

with $\varepsilon > 0$. Now we can show that if we have

(6.23) $$\|F_\mathbb{S} - \Phi\| \ll \beta + \gamma^\vartheta,$$

then we must have $\vartheta \leq 1/2$.

We prove the above claim by contradiction. Namely, assume that the contrary holds, that is, that (6.23) holds with some $\vartheta > 1/2$. Clearly, $\beta = \mathbb{E}|X_1|^3/\sqrt{n} \ll n^{-1/2}$, $\gamma = \operatorname{var} \mathbb{T} = O(n^{-1})$, $\kappa_1 = 0$, $\kappa_2 = \varepsilon n^{-1/2}$, $G_1(x) = 0$, and $G_1(x) = \varepsilon n^{-1/2} \Phi'''(x)$. In view of Theorem 1 and the assumption (6.23), we have

$$|G_2(x)| = \|G_2\| \leq \|F_S - \Phi\| + \|F_S - \Phi - G_2\| \ll \beta + \gamma^\vartheta + \gamma \ln^2(\gamma).$$

Multiplying $\sqrt{n}$ on both sides, we get

$$|\varepsilon \Phi'''(x)| \ll 1 + \sqrt{n}\gamma^\vartheta + \sqrt{n}\gamma \ln^2(\gamma).$$

Letting $n \to \infty$, and in view of $\gamma = O(n^{-1})$ and $\vartheta > 1/2$, we get $\varepsilon|\Phi''(x)| \ll 1$ which contradicts the assumption that $\varepsilon > 0$ is an arbitrary positive number. Thus we must have $\vartheta \leq 1/2$.

6.4. *Proof of Theorem 5.* Let $\eta$ be a r.v. from $N(0,1)$. Note that the left-hand side of the inequality in Theorem 5 can be written as

$$\mathbb{E}f(\mathbb{S}) = \int_\mathbb{R} f(u) \, d\mathbb{P}\{\mathbb{S} \leq u\},$$

$$\mathbb{E}f(\eta) + \kappa_1 \mathbb{E}f''(\eta) + \cdots + \kappa_k \mathbb{E}f^{(k+1)}(\eta) = \int_\mathbb{R} f(u) \, dN_k(u).$$

Expanding $f(\mathbb{S}) = f(\mathbb{L} + \mathbb{T})$ in powers of $\mathbb{T}$, we get

$$\mathbb{E}f(\mathbb{S}) = \mathbb{E}f(\mathbb{L}) + \mathbb{E}f'(\mathbb{L})\mathbb{T} + \vartheta, \qquad \vartheta \stackrel{\text{def}}{=} \mathbb{E}(1-\tau)f''(\mathbb{L} + \tau\mathbb{T})\mathbb{T}^2,$$

where $\tau$ is uniformly distributed on $[0,1]$, and is independent of all other r.v.'s. To prove the theorem it suffices to check that

(6.24) $$\|\mathbb{E}f(\mathbb{L}) - \mathbb{E}f(\eta)\| \ll \|f'''\|\beta,$$

(6.25) $$\|\mathbb{E}\mathbb{T}f'(\mathbb{L}) - I\| \ll C(\beta + \gamma),$$

(6.26) $$\|\vartheta\| \leq \|f''\|\gamma,$$

where

(6.27) $$I = \kappa_1 \mathbb{E}f''(\eta) + \cdots + \kappa_k \mathbb{E}f^{(k+1)}(\eta) \stackrel{\text{def}}{=} I_1 + \cdots + I_k.$$

The estimate (6.24) is well known since $\mathbb{L} = L_1 + \cdots + L_n$ is a sum of i.i.d. r.v.'s (see, e.g., [5]).

The bound (6.26) is obvious. Indeed, we have $\|\vartheta\| \leq \|f''\|\mathbb{E}\mathbb{T}^2 = \|f''\|\gamma$.

It remains to prove (6.25). Using $\mathbb{T} = \mathbb{T}_1 + \cdots + \mathbb{T}_k$, the linear structure of $\mathbb{T}_s$, and the i.i.d. assumption, we can write

(6.28) $$\mathbb{E}\mathbb{T}f'(\mathbb{L}) = J_1 + \cdots + J_k, \qquad J_s \stackrel{\text{def}}{=} \binom{n}{s} \mathbb{E}T_{1\cdots s}f'(\mathbb{L}).$$

In view of (6.27) and (6.28), the proof of (6.25) reduces to checking that

(6.29) $$\|I_s - J_s\| \ll C(\beta + \gamma).$$

Let us split $\mathbb{L} = W + R$, where $W = L_1 + \cdots + L_s$ and $R = L_{s+1} + \cdots + L_n$. Writing $K_s = \kappa_s \mathbb{E} f^{(s+1)}(R)$, instead of (6.29), it suffices to prove that

(6.30) $$\|I_s - K_s\| \ll C(\beta + \gamma)$$

and

(6.31) $$\|K_s - J_s\| \ll C(\beta + \gamma).$$

PROOF OF (6.30). Note that $I_s - K_s = \kappa_s(\mathbb{E} f^{(s+1)}(\eta) - \mathbb{E} f^{(s+1)}(R))$. Since $\mathbb{L} = R + W$, to prove (6.30) it suffices to show that

(6.32) $$|\kappa_s| \|\mathbb{E} f^{(s+1)}(\eta) - \mathbb{E} f^{(s+1)}(\mathbb{L})\| \ll C(\beta + \gamma)$$

and

(6.33) $$|\kappa_s| \|\mathbb{E} f^{(s+1)}(R + W) - \mathbb{E} f^{(s+1)}(R)\| \ll C(\beta + \gamma).$$

Let us prove (6.32). We consider the cases $\beta \geq 1$ and $\beta < 1$ separately. Using $\kappa_s^2 \leq \gamma$ (see Lemma 1), in the case $\beta \geq 1$ we have

$$|\kappa_s| \|\mathbb{E} f^{(s+1)}(\eta) - \mathbb{E} f^{(s+1)}(\mathbb{L})\| \leq 2\|f^{(s+1)}\|\sqrt{\gamma} \ll C\sqrt{\gamma\beta} \leq C(\gamma + \beta).$$

In the case $\beta \leq 1$ we use the bound (6.24) replacing $f$ by $f^{(s+1)}$. Since $\kappa_s^2 \leq \gamma$, we get

$$|\kappa_s| \|\mathbb{E} f^{(s+1)}(\eta) - \mathbb{E} f^{(s+1)}(\mathbb{L})\| \ll \|f^{(s+4)}\|\sqrt{\gamma}\beta \leq C\sqrt{\gamma\beta} \leq C(\gamma + \beta).$$

Let us prove (6.33). We apply (6.1) replacing $g$ by $f^{(s+1)}$, the variable $x$ by $R$, and $h$ by $W$, respectively. Since $\kappa_s^2 \leq \gamma$, we have

(6.34) $$|\kappa_s| \|\mathbb{E} f^{(s+1)}(R + W) - \mathbb{E} f^{(s+1)}(R)\|$$
$$\ll (\|f^{(s+1)}\| + \|f^{(s+2)}\|)\gamma^{1/2} \mathbb{E}|W|^{1/2}.$$

Using $|W|^{1/2} \leq |L_1|^{1/2} + \cdots + |L_s|^{1/2}$ and $\mathbb{E} L_m^2 = 1/n$, we derive $\mathbb{E}|W|^{1/2} \ll n^{-1/4} \leq \beta^{1/2}$. In view of (6.34), an application of $\sqrt{\gamma\beta} \leq \gamma + \beta$ yields (6.33).

PROOF OF (6.31). Using short Taylor expansions, we can represent $J_s$ as

(6.35) $$J_s = \binom{n}{s} \mathbb{E} L_1 \cdots L_s T_{1 \ldots s} f^{(s+1)}(V + R)$$

with $V = \tau_1 L_1 + \cdots + \tau_s L_s$. Representation (6.35) can be proven in $s$ steps. Let us consider details related only to the first step. The degeneracy property of the kernels implies that

(6.36) $$\mathbb{E} T_{1 \ldots s} f'(L_2 + \cdots + L_n) = 0$$

since the conditional expectation $\mathbb{E}(T_{1\cdots s}|X_2,\ldots,X_{s+1}) = 0$, and $L_2 + \cdots + L_n$ is independent of $X_1$. Using (6.36), expanding in powers of $L_1$ we get

$$J_s = \binom{n}{s} \mathbb{E} L_1 T_{1\cdots s} f''(\tau_1 L_1 + L_2 + \cdots + L_n).$$

Proceeding in a similar way with $L_2,\ldots,L_s$ in place of $L_1$, we arrive at (6.35).

Using the definition of $\kappa_s$, and the expression (6.35) for $J_s$, we can write $K_s - J_s = \binom{n}{s} \mathbb{E} L_1 \cdots L_s T_{1\cdots s}(f^{(s+1)}(R) - f^{(s+1)}(V+R))$. By an application of (6.1) with $g = f^{(s+1)}$, $x = R$ and $h = V$, it follows that

$$(6.37) \qquad \|K_s - J_s\| \leq C \binom{n}{s} \mathbb{E}|L_1 \cdots L_s T_{1\cdots s}||V|^{1/2}.$$

Using the inequality $|V|^{1/2} \leq |L_1|^{1/2} + \cdots + |L_s|^{1/2}$, the Hölder inequality, the assumption $\mathbb{E}|L_1|^3 = \beta/n$ and $\mathbb{E} L_1^2 = 1/n$, we get

$$\mathbb{E}|L_1 \cdots L_s T_{1\cdots s}||V|^{1/2} \leq s\mathbb{E}|L_1|^{3/2}|L_2 \cdots L_s T_{1\cdots s}|$$
$$(6.38) \qquad \ll (\mathbb{E}|L_1|^3 L_2^2 \cdots L_s^2 \mathbb{E} T_{1\cdots s}^2)^{1/2}$$
$$= n^{-s/2}(\beta \mathbb{E} T_{1\cdots s}^2)^{1/2}.$$

Using $\mathbb{E} T_{1\cdots s}^2 = \mathbb{E} \mathbb{T}_s^2 / \binom{n}{s} \leq \gamma / \binom{n}{s}$, relations (6.37) and (6.38) yield $\|K_s - J_s\| \ll C\sqrt{\alpha}\sqrt{\beta\gamma}$, $\alpha = \binom{n}{s}/n^s$. Noting that $\alpha \leq 1$ and $\sqrt{\beta\gamma} \leq \beta + \gamma$, we derive $\|K_s - J_s\| \ll C(\beta + \gamma)$ which concludes the proof of (6.31).

6.5. *Proof of Theorem 6.* Let $\eta$ be a r.v. from $N(0,1)$. Note that the left-hand side of the inequality in Theorem 6 can be written as

$$\mathbb{E} f(\mathbb{S}) = \int_{\mathbb{R}} f(u)\, d\mathbb{P}\{\mathbb{S} \leq u\},$$
$$\mathbb{E} f(\eta) + \kappa_1 \mathbb{E} f''(\eta) + \cdots + \kappa_p \mathbb{E} f^{(p+1)}(\eta)$$
$$= \int_{\mathbb{R}} f(u)\, dN_p(u).$$

We consider the cases $\alpha = 1$ and $1 < \alpha < 2$, separately.

First consider the case $\alpha = 1$. In this case $I^{(\alpha)} = 0$. Similar to (6.1), one can check that

$$(6.39) \qquad f(\mathbb{S}) = f(\mathbb{L}) + \vartheta \qquad \text{with } \|\vartheta\| \ll \|f\|^{1-\alpha} \|f'\|^\alpha |\mathbb{T}|^\alpha.$$

It is easy to show that

$$\mathbb{E}|\mathbb{T}_p|^\alpha \ll \gamma_p^{(\alpha)} = \binom{n}{p} \mathbb{E}|T_{1\cdots p}|^\alpha,$$
$$(6.40)$$
$$\mathbb{E}|\mathbb{T}|^\alpha \ll \gamma^{(\alpha)} = \gamma_1^{(\alpha)} + \cdots + \gamma_k^{(\alpha)}.$$

Then the proof follows from (6.24), (6.39) and (6.40).

Now consider the case $1 < \alpha < 2$. Similar to (6.1), we have

$$f(\mathbb{S}) = f(\mathbb{L}) + f'(\mathbb{L})\mathbb{T} + \vartheta \qquad \text{with } \|\vartheta\| \ll (\|f'\| + \|f''\|)|\mathbb{T}|^\alpha.$$

From (6.40), we get $\mathbb{E}|\mathbb{T}|^\alpha \ll \gamma^{(\alpha)}$. Also from (6.24) it follows easily that $\|\mathbb{E}f(\mathbb{L}) - \mathbb{E}f(\eta)\| \ll \|f'''\|\beta$. It remains to show that

$$\|\mathbb{E}\mathbb{T}f'(\mathbb{L}) - I^{(\alpha)}\| \ll C(\beta + \gamma^{(\alpha)}), \tag{6.41}$$

where $I^{(\alpha)} = \kappa_1 \mathbb{E}f''(\eta) + \cdots + \kappa_p \mathbb{E}f^{(p+1)}(\eta)$. We have

$$\mathbb{E}\mathbb{T}f'(\mathbb{L}) = J_1 + \cdots + J_k, \ J_s \stackrel{\text{def}}{=} \binom{n}{s} \mathbb{E}T_{1\cdots s}f'(\mathbb{L}). \tag{6.42}$$

In view of (6.41) and (6.42), the proof of the theorem reduces to checking

$$\|I_s - J_s\| \ll C(\beta + \gamma^{(\alpha)}) \qquad \text{for } 1 \leq s \leq p, \tag{6.43}$$

with $I_s = \kappa_s \mathbb{E}f^{(s+1)}(\eta)$, and

$$\|J_s\| \ll C(\beta + \gamma^{(\alpha)}) \qquad \text{for } p < s \leq k, \tag{6.44}$$

where $p$ is an integer satisfying the condition of the theorem.

*Estimation of $J_s$ for $s > p$.* We have to prove (6.44). Consider the difference operator $\delta_m$ such that $\delta_m f(x) = f(x + L_m) - f(x)$. Arguments similar to those used to derive (6.35) allow us to write

$$J_s = \binom{n}{s} \mathbb{E}T_{1\cdots s}\delta_1 \cdots \delta_s f'(R) \tag{6.45}$$

with $R = L_{s+1} + \cdots + L_n$. Similar to (6.1), we have

$$\|\delta_1 \cdots \delta_s f'(R)\| \ll \|f'\| + \|f^{(s+1)}\||L_1|^\varepsilon \cdots |L_s|^\varepsilon, \qquad 0 \leq \varepsilon \leq 1.$$

From this and the Hölder inequality with exponents $q = \alpha/(\alpha - 1)$ and $\alpha$, so that $1/q + 1/\alpha = 1$, we have

$$\|J_s\| \ll C\binom{n}{s} \mathbb{E}|L_1|^\varepsilon \cdots |L_s|^\varepsilon |T_{1\cdots s}|$$

$$\ll C\binom{n}{s} (\mathbb{E}|L_1|^{\varepsilon q} \cdots |L_s|^{\varepsilon q})^{1/q} (\mathbb{E}|T_{1\cdots s}|^\alpha)^{1/\alpha}$$

$$= C\left(\binom{n}{s}(\mathbb{E}|L_1|^{\varepsilon q})^s\right)^{1/q} (\gamma^{(\alpha)})^{1/\alpha}$$

$$\ll C\binom{n}{s}(\mathbb{E}|L_1|^{\varepsilon q})^s + C\gamma^{(\alpha)}.$$

To complete the proof of (6.44) we need to prove $\binom{n}{s}(\mathbb{E}|L_1|^{\varepsilon q})^s \ll \beta$. Choose $\varepsilon = (2 + 1/s)/q$. The condition $s > p$ guaranties that $\varepsilon \leq 1$. Using $n\mathbb{E}|L_1|^{2+1/s} \leq$ $\beta^{1/s}$ [see (b) of Lemma 1], we have $\binom{n}{s}(\mathbb{E}|L_1|^{\varepsilon q})^s \leq (n\mathbb{E}|L_1|^{2+1/s})^s \leq \beta$, which completes the proof for $s > p$.

*Estimation of $J_s$ for $1 \leq s \leq p$.* We have to prove (6.43). There is something to prove only if $3/2 < \alpha < 2$. Indeed, in the case $1 < \alpha \leq 3/2$, the integer $p$ from the condition of the theorem satisfies $p \leq 0$, and we have no $J_s$ to estimate.

Let us split $\mathbb{L} = W + R$ where $W = L_1 + \cdots + L_s$ and $R = L_{s+1} + \cdots + L_n$. Writing $K_s = \kappa_s \mathbb{E} f^{(s+1)}(R)$ instead of (6.43) it suffices to prove that $\|I_s - K_s\|$ and $\|K_s - J_s\|$ are bounded from above like $\|I_s - J_s\|$ in (6.43).

*Estimation of $\|I_s - K_s\|$.* Note that $I_s - K_s = \kappa_s \mathbb{E}(f^{(s+1)}(\eta) - f^{(s+1)}(R))$. Since $\mathbb{L} = R + W$, it suffices to show that

$$(6.46) \qquad |\kappa_s| \|f^{(s+1)}(\eta) - f^{(s+1)}(\mathbb{L})\| \ll C(\beta + \gamma^{(\alpha)})$$

and

$$(6.47) \qquad |\kappa_s| \|\mathbb{E} f^{(s+1)}(R + W) - \mathbb{E} f^{(s+1)}(R)\| \ll C(\beta + \gamma^{(\alpha)}).$$

Let us prove (6.46). We consider the cases $\beta \geq 1$ and $\beta \leq 1$ separately. Using $|\kappa_s| \ll \beta^{\delta s} + \gamma^{(\alpha)}$ and $\delta s \leq 1$ [see (d) of Lemma 1], in the case of $\beta \geq 1$, we have

$$|\kappa_s| \|\mathbb{E} f^{(s+1)}(\eta) - \mathbb{E} f^{(s+1)}(\mathbb{L})\| \leq 2\|f^{(s+1)}\|(\beta^{\delta s} + \gamma^{(\alpha)}) \ll C(\beta + \gamma^{(\alpha)}).$$

In the case of $\beta \leq 1$, we combine (d) of Lemma 1 with the bound (6.24) replacing $f$ by $f^{(s+1)}$. We get

$$|\kappa_s| \|\mathbb{E} f^{(s+1)}(\eta) - \mathbb{E} f^{(s+1)}(\mathbb{L})\| \ll \|f^{(s+4)}\|(\beta^{\delta s} + \gamma^{(\alpha)})\beta \leq C(\beta + \gamma^{(\alpha)}).$$

Let us now prove (6.47). Again we consider the cases of $\beta \geq 1$ and $\beta \leq 1$ separately. We omit details in the case of $\beta \geq 1$ since they are very similar to those in the proof of (6.46). In the case of $\beta \leq 1$, we use (d) of Lemma 1 and apply $\|f^{(s+1)}(R + W) - \mathbb{E} f^{(s+1)}(R)\| \leq C\|W\|$ to get

$$(6.48) \quad |\kappa_s| \|\mathbb{E} f^{(s+1)}(R + W) - \mathbb{E} f^{(s+1)}(R)\| \ll C(\beta^{\delta s} + \gamma^{(\alpha)}) \mathbb{E}|W|.$$

Using $|W| \leq |L_1| + \cdots + |L_s|$ and $\mathbb{E} L_m^2 = 1/n$, we derive $\mathbb{E}|W| \ll n^{-1/2} \leq \beta$. In view of our assumption $\beta \leq 1$, the inequality (6.48) yields (6.47).

*Estimation of $\|K_s - J_s\|$.* By (6.35) we have

$$J_s = \binom{n}{s} \mathbb{E} L_1 \cdots L_s T_{1 \cdots s} f^{(s+1)}(V + R)$$

with $V = \tau_1 L_1 + \cdots + \tau_s L_s$. This representation, combined with the definition $\kappa_s$, allows one to write

$$K_s - J_s = \binom{n}{s} \mathbb{E} L_1 \cdots L_s T_{1 \cdots s}(f^{(s+1)}(R) - f^{(s+1)}(V + R)).$$

Again we consider the cases of $\beta \geq 1$ and $\beta \leq 1$ separately. We omit details in the case of $\beta \geq 1$ since they are very similar to those in the proof of (6.46). In the case of $\beta \leq 1$, we apply

$$\|\mathbb{E}f^{(s+1)}(R) - f^{(s+1)}(R+V)\| \leq C\|V\|^\varepsilon, \qquad 0 \leq \varepsilon \leq 1.$$

From the inequality $|V|^\varepsilon \ll |L_1|^\varepsilon + \cdots + |L_s|^\varepsilon$, the i.i.d. assumption, the Hölder inequality with exponents $q = \alpha/(\alpha - 1)$ and $\alpha$, and also $\binom{n}{s} \leq n^s$, we have

$$\begin{aligned}
\|K_s - J_s\| &\ll C\binom{n}{s}\mathbb{E}|L_1 \cdots L_s T_{1\cdots s}||V|^\varepsilon \\
&\leq Cn^s s\mathbb{E}|L_1|^{1+\varepsilon}|L_2 \cdots L_s T_{1\cdots s}| \\
&\ll Cn^s(\mathbb{E}|L_1|^{(1+\varepsilon)q}(\mathbb{E}|L_1|^q)^{s-1})^{1/q}(\mathbb{E}|T_{1\cdots s}|^\alpha)^{1/\alpha} \\
&\ll Cn^s\mathbb{E}|L_1|^{(1+\varepsilon)q}(\mathbb{E}|L_1|^q)^{s-1} + Cn^s\mathbb{E}|T_{1\cdots s}|^\alpha \\
&= Cn^s\mathbb{E}|L_1|^{(1+\varepsilon)q}(\mathbb{E}|L_1|^q)^{s-1} + C\gamma^{(\alpha)}.
\end{aligned}$$

We choose $\varepsilon$ so that $(1+\varepsilon)q = 3$. In view of $3/2 < \alpha < 2$, the number $\varepsilon$ satisfies $0 < \varepsilon < 1$. Since $n\mathbb{E}|L_1|^3 = \beta$, to conclude the estimation of $\|K_s - J_s\|$ it suffices to verify the inequality $n\mathbb{E}|L_1|^q \leq 1$. Recalling that $\beta \leq 1$ and applying Lemma 1, we have $n\mathbb{E}|L_1|^q \leq \beta^{q-2} \leq 1$, which completes the estimation of $\|K_s - J_s\|$.

**Acknowledgments.** The authors would like to thank the Editor, Associate Editor and referees for their comments and criticisms, which led to this much improved version.

[8] FELLER, W. (1971). *An Introduction to Probability Theory and Its Applications. Vol. II*, 2nd ed. Wiley, New York. MR0270403
[9] FRIEDRICH, K. O. (1989). A Berry–Esseen bound for functions of independent random variables. *Ann. Statist.* **17** 170–183. MR981443
[10] GINÉ, E., LATAŁA, R. and ZINN, J. (2000). Exponential and moment inequalities for $U$-statistics. In *High Dimensional Probability, II (Seattle, WA, 1999). Progress in Probability* **47** 13–38. Birkhäuser, Boston, MA. MR1857312
[11] HOEFFDING, W. (1948). A class of statistics with asymptotically normal distribution. *Ann. Math. Statist.* **19** 293–325. MR0026294
[12] JING, B.-Y. and WANG, Q. (2003). Edgeworth expansion for $U$-statistics under minimal conditions. *Ann. Statist.* **31** 1376–1391. MR2001653
[13] KARLIN, S. and RINOTT, Y. (1982). Applications of ANOVA type decompositions for comparisons of conditional variance statistics including jackknife estimates. *Ann. Statist.* **10** 485–501. MR653524
[14] KOROLYUK, V. S. and BOROVSKIKH, YU. V. (1985). Approximation of nondegenerate $U$-statistics. *Theory Probab. Appl.* **30** 439–450.
[15] KOROLJUK, V. S. and BOROVSKICH, Y. V. (1994). *Theory of $U$-Statistics. Mathematics and Its Applications* **273**. Kluwer, Dordrecht. MR1472486
[16] PECCATI, G. (2004). Hoeffding-ANOVA decompositions for symmetric statistics of exchangeable observations. *Ann. Probab.* **32** 1796–1829. MR2073178
[17] VAN ZWET, W. R. (1984). A Berry–Esseen bound for symmetric statistics. *Z. Wahrsch. Verw. Gebiete* **66** 425–440. MR751580

V. BENTKUS
INSTITUTE OF MATHEMATICS
 AND INFORMATICS
AKADEMIJOS 4, VILNIUS
LITHUANIA
E-MAIL: bentkus@ktl.mii.lt

B.-Y. JING
DEPARTMENT OF MATHEMATICS
HONG KONG UNIVERSITY OF SCIENCE
 AND TECHNOLOGY
CLEAR WATER BAY
KOWLON
HONG KONG
E-MAIL: majing@ust.hk

W. ZHOU
DEPARTMENT OF STATISTICS
 AND APPLIED PROBABILITY
NATIONAL UNIVERSITY OF SINGAPORE
6 SCIENCE DRIVE 2
SINGAPORE 117546
SINGAPORE
E-MAIL: stazw@nus.edu.sg